\documentclass{amsart}

\usepackage[english]{babel}
\usepackage{amsmath,amssymb,amsfonts,amsthm}

\usepackage[centering]{geometry}
\usepackage{color}
\usepackage[utf8]{inputenc}

\usepackage{hyperref}
\hypersetup{colorlinks   = true,
  		    citecolor    = green,
   		    linkcolor    = blue}
\usepackage{cleveref}
\usepackage{mathtools}
\usepackage{ mathrsfs }
\usepackage{environ}
\usepackage{caption}
\usepackage{systeme}
\usepackage{subcaption}
\usepackage{listings}
\newcommand{\R}{\ensuremath{\mathbb{R}}}

\def\<{\left < }
\def\>{\right >}
\def\({\left ( }
\def\){\right )}

\newcommand{\de}{\mathrm{d}}

\newtheorem{theorem}{Theorem}[section]   %% definition of theorem environment
\newtheorem*{theorem*}{Theorem}          %% a theorem environment without numbering

\newtheorem{proposition}[theorem]{Proposition}

\theoremstyle{definition}
\newtheorem{definition}[theorem]{Definition}
\newtheorem{corollary}[theorem]{Corollary}

\newtheorem*{remark}{Remark}

\title[A geometric study of marginally trapped surfaces]{A geometric study of marginally trapped surfaces in space forms and Robertson-Walker spacetimes -- an overview}

\author{Kristof Dekimpe}
\author{Joeri Van der Veken}
\address{KU Leuven, Department of Mathematics, Celestijnenlaan 200B - Box 2400, 3001 Leuven, Belgium}
\email{kristof.dekimpe@kuleuven.be}
\email{joeri.vanderveken@kuleuven.be}
\thanks{The second author is supported by the Excellence Of Science project G0H4518N of the Belgian government and by project 3E160361 of the KU Leuven Research Fund. Both authors are supported by the collaboration project G0F2319N of the Research Foundation -- Flanders (FWO) and the National Natural Science Foundation of China (NSFC)}


\begin{document}

\maketitle

\begin{abstract}
A marginally trapped surface in a spacetime is a Riemannian surface whose mean curvature vector is lightlike at every point. In this paper we give an up-to-date overview of the differential geometric study of these surfaces in Minkowski, de Sitter, anti-de Sitter and Robertson-Walker spacetimes. We give the general local descriptions proven by Anciaux and his coworkers as well as the known classifications of marginally trapped surfaces satisfying one of the following additional geometric conditions: having positive relative nullity, having parallel mean curvature vector field, having finite type Gauss map, being invariant under a one-parameter group of ambient isometries, being isotropic, being pseudo-umbilical. Finally, we provide examples of constant Gaussian curvature marginally trapped surfaces and state some open questions. 
\end{abstract}

%%%%%%%%%%%%%%%%%%%%%%%%%%%%%%%%%%%%%%%%%%%%%%%%%%%%%%%%%%%

\section{Introduction}

Trapped surfaces were introduced by Sir Roger Penrose in \cite{Penrose} %%MDPI: ref citations are not permitted in Abstract, please move [1,2] into main text 
and play an important role in cosmology. From a purely differential geometric point of view, a marginally trapped surface in a spacetime is a Riemannian surface whose mean curvature vector field is lightlike at every point, i.e., for every point $p$ of the surface, the mean curvature vector $H(p)$ satisfies $\langle H(p),H(p) \rangle = 0$ and $H(p) \neq 0$, while every non-zero vector $v$ tangent to the surface satisfies $\langle v,v \rangle > 0$.

Since 2007, several authors have studied marginally trapped surfaces in spacetimes and their generalizations to higher signatures, dimensions and codimensions from a geometric point of view. Most of the results give a complete classification of marginally trapped surfaces in a specific spacetime under one or more additional geometric conditions, such as having positive relative nullity \cite{Positiverelativenullity}, having parallel mean curvature vector field \cite{Parallelmeancurvature}, having finite type Gauss map \cite{Gaussmap}, being invariant under certain $1$-parameter groups of isometries \cite{Boost, Rotation, Screw} or being isotropic \cite{Isotropy}. Several of these results and the above mentioned generalizations are due to Bang-Yen Chen and his collaborators and we should also mention his 2009 overview paper on the topic \cite{Chen2009}. In 2015, Henri Anciaux and his collaborators gave a local description of any marginally trapped surface (and even codimension two submanifold) in a Lorentzian space form or a Robertson-Walker spacetime in \cite{localRobertson, localspace}, without requiring any additional properties. While this is the most general result, and in some sense the best one can hope for in this context, the previously mentioned results are still of great value, since they often give more explicit descriptions under the additional condition at hand. Also, it is not always easy to find the surfaces with a certain property from the general description, such as those with constant Gaussian curvature.

In this paper we give an overview of all the above mentioned results. We work in reverse chronological order, in the sense that, after the preliminaries (Section \ref{sec:preliminaries}), we state the general results by Anciaux et al. (Section \ref{sec:local_description}). After that, we discuss the classifications under additonal assumptions which can be found in the literature (Section \ref{sec:posrelnul}, \ref{sec:pmc}, \ref{sec:invariant}, \ref{sec:isotropic}). In Section \ref{sec:cgc}, we discuss marginally trapped surfaces with constant Gaussian curvature. The result is limited to finding the examples in previously obtained lists and hence does not provide a complete classification. We end the paper with some open problems (Section \ref{sec:openproblems}).

%%%%%%%%%%%%%%%%%%%%%%%%%%%%%%%%%%%%%%%%%%%%%%%%%%%%%%%%%%%

\section{Preliminaries}\label{sec:preliminaries}

A pseudo-Riemannian manifold is a manifold with a metric which is not required to be positive definite, but merely non-degenerate. Consequently, one can distinguish three types of tangent vectors to such a manifold: a vector $v$ is \textit{spacelike} if $\langle v,v \rangle > 0$ or $v=0$, \textit{timelike} if $\langle v,v \rangle < 0$ and \textit{lightlike} or \textit{null} if $\langle v,v \rangle = 0$ but $v \neq 0$.

Let $S$ be a Riemannian or spacelike $n$-dimensional submanifold of a pseudo-Riemannian manifold $M$ and denote by $\nabla$ and $\widetilde{\nabla}$ the Levi-Civita connections of $S$ and $M$ respectively. The formulas of Gauss and Weingarten state respectively that
\begin{align*}
& \widetilde{\nabla}_X Y=\nabla_X Y + h(X,Y), \\
& \widetilde{\nabla}_X \xi = -A_{\xi}X + \nabla^{\perp}_X \xi
\end{align*}
for vector fields $X$ and $Y$ tangent to $S$ and a vector field $\xi$ normal to $S$. Here, $h$ is a symmetric $(1,2)$-tensor field taking values in the normal bundle, called the \textit{second fundamental form}, $A_{\xi}$ is a symmetric $(1,1)$-tensor field, called the \textit{shape operator} associated to $\xi$, and $\nabla^{\perp}$ is the \textit{normal connection}. The \textit{mean curvature vector} at a point $p\in S$  is defined by
\begin{align*}
H(p) = \frac{1}{n} \sum_{i=1}^n h(e_i, e_i),
\end{align*}
where $\left\{ e_1,\ldots,e_n \right\}$ is an orthonormal basis of $T_pS$. We say that $S$ is \textit{marginally trapped} if its mean curvature vector $H(p)$ is null for every point $p$ of $S$ and we say that $S$ has \textit{null second fundamental form} if $h(X,Y)$ is a null vector for all tangent vectors $X$ and $Y$ to $S$. It is clear that submanifolds with null second fundamental form are marginally trapped, but the converse is not necessarily true. Submanifolds for which the mean curvature vanishes at some points are not marginally trapped, since the zero vector is not a null vector. However, since we stay close to the original sources, some of the classifications in this paper include conditions preventing the vanishing of the mean curvature vector field, while some others don't (see also the remark after Theorem \ref{T:LocalMinkowski}). For an explicit example of a submanifold, it is not hard to check whether the mean curvature vanishes at some points or not.

We consider two types of ambient Lorentzian manifolds in this paper. The first type are those with the highest degree of symmetry, namely Lorentzian real space forms. We recall the definition of a (pseudo-)Riemannian space form of any index. Let $\mathbb{R}^n_s$ denote $\R^n = \{(x_1,\ldots,x_n) \ | \ x_1,\ldots,x_n \in~\R \}$ equipped with the inner product
\begin{align*}
\langle (x_1,\ldots,x_n) , (y_1,\ldots,y_n) \rangle = -x_1 y_1 - \ldots - x_s y_s + x_{s+1} y_{s+1} + \ldots + x_n y_n.
\end{align*}
For $s=0$, the space $\R^n_0 = \R^n$ is just the Euclidean space of dimension $n$ and for $s>1$, we call $\R^n_s$ the \textit{pseudo-Euclidean space} of dimension $n$ and index $s$. It is a flat manifold, i.e., a pseudo-Riemannian manifold with constant sectional curvature $0$. We now define
\begin{align*}
& S^n_s(x_0,c) =\left\{x\in\mathbb{R}^{n+1}_s \ | \ \left<x-x_0,x-x_0 \right> = 1/c \right\} \mbox{ for } x_0\in \mathbb{R}^{n+1}_s \mbox{ and }c>0,\\
& H^n_s(x_0,c) =\left\{x\in\mathbb{R}^{n+1}_{s+1} \ | \ \left<x-x_0,x-x_0 \right> = 1/c \right\} \mbox{ for } x_0\in \mathbb{R}^{n+1}_{s+1} \mbox{ and }c<0.
\end{align*}
The manifolds $S^n_s(x_0,c)$ and $H^n_s(x_0, c)$, equipped with the induced metrics from $\R^{n+1}_s$ and $\R^{n+1}_{s+1}$ respectively, are complete pseudo-Riemannian manifolds with constant sectional curvature $c$. We refer to $\R^n_s$, $S^n_s(x_0,c)$ and $H^n_s(x_0,c)$ as \textit{real space forms} of dimension $n$ and index $s$. We simply denote $S^n_s(x_0,c)$ and $H^n_s(x_0,c)$ by $S^n_s(c)$ and $H^n_s(c)$ when $x_0$ is the origin. In the rest of this paper, we will often use the unified notation $Q^n_s(c)$ to denote the following:
\begin{align*}
Q^n_s(c) = \left\{ \begin{array}{ll}
S^n_s(c) & \mbox{if $c>0$}, \\
\R^n_s & \mbox{if $c=0$}, \\
H^n_s(c) & \mbox{if $c<0$}. \\
\end{array}\right.
\end{align*}
If the index $s=0$, we denote the Riemannian real space form $Q^n_0(c)$ by $Q^n(c)$ and if the index $s=1$, we call $Q^n_1(c)$ a \textit{Lorentzian real space form}. In particular, the four-dimensional Lorentzian real space forms  $\mathbb{R}^4_1$, $S^4_1(1)$ and $H^4_1(-1)$ are known as the \textit{Minkowksi spacetime}, the \textit{de Sitter spacetime} and the \textit{anti-de Sitter spacetime}.

The second type of ambient Lorentzian manifolds that we consider are warped products
\begin{align*}
L^n_1(f,c)=I\times_f Q^{n-1}(c), 
\end{align*}
where $I\subset\mathbb{R}$ is an open interval, $f: I \to \R$ is a smooth positive function and $Q^{n-1}(c)$ is a Riemannian real space form with constant curvature $c\in\left\{-1,0,1 \right\}$. The metric of $L^n_1(f,c)$ is given by
\begin{align*}
\langle \,\cdot\,,\cdot \,\rangle = -\de t^2+f^2(t) \langle \,\cdot\,,\cdot \,\rangle_c,
\end{align*}
where $t$ is a coordinate on $I$ and $\langle \,\cdot\,,\cdot \,\rangle_c$ is the metric of $Q^{n-1}(c)$. If the warping function $f$ is constant, then $L^n_1(f,c)$ is the Lorentzian product of $(I,-\de t^2)$ and a Riemannian real space form. For general $f$ but $n=4$, the Lorentzian manifolds $L^4_1(f,c)$ are known as \textit{Robertson-Walker spacetimes}. The next remark, which can for example be found in \cite{Robertson} in a slightly different form, determines when a Robertson-Walker spacetime has constant sectional curvature. 
\begin{remark}
A Robertson-Walker spacetime $L^4_1(f,c)$ has constant sectional curvature $K$ if and only if the warping function $f$ satisfies
\begin{align*}
f^2 K=ff''=\left(f' \right)^2+c.
\end{align*}
We thus have the following:
\begin{itemize}
\item[(1)]$L^4_1(f,c)$ is flat if and only if $f(t)=at+b$, with $a^2=-c$;
\vskip.05in
\item[(2)]$L^4_1(f,c)$ has constant sectional curvature $K>0$ if and only if 
\begin{align*}
f(t)=a\cosh(\sqrt{K}t)+b\sinh(\sqrt{K}t), \mbox{ with } a^2-b^2=\frac{c}{K}\mbox{;}
%f(t)=ae^{\sqrt{K}t}+be^{-\sqrt{K}t} \mbox{ , with }ab=\frac{c}{4K};
\end{align*}
\item[(3)]$L^4_1(f,c)$ has constant sectional curvature $K<0$ if and only if 
\begin{align*}
f(t)=a\cos(\sqrt{-K}t)+b\sin(\sqrt{-K}t), \mbox{ with } a^2+b^2=\frac{c}{K}.
\end{align*}
\end{itemize}
\end{remark}
\noindent We end this section by defining the following family of functions on a space of type $L^n_1(f,c)$ which will appear in some of the results:
\begin{equation} \label{eq:theta}
\theta : L^n_1(f,c) \to \R : t \mapsto \int_{t_0}^{t}\frac{\de s}{f(s)}
\end{equation}
for $t_0 \in I$.

%%%%%%%%%%%%%%%%%%%%%%%%%%%%%%%%%%%%%%%%%%%%%%%%%%%%%%%%%%%

\section{Local Description of Codimension Two Marginally Trapped Submanifolds}\label{sec:local_description}

In this section, we give a local description of \textit{any} codimension two marginally trapped submanifold of a Lorentzian space form or a space of type $L^n_1(f,c)$. In particular, these descriptions hold for marginally trapped surfaces in Minkowski spacetime, de Sitter spacetime, anti-de Sitter spacetime and Robertson-Walker spacetimes. The results were proven in \cite{localRobertson,localspace}.

The results use the notion of the Gauss map of a hypersurface of a Riemannian space form.

\begin{definition}
Let $S$ be an immersed hypersurface of a Riemannian real space from $Q^{n+1}(c)$ and denote by $\nu$ a unit normal vector field along the immersion. Such a vector field always exists locally and it exists globally if $S$ is orientable. Since $Q^{n+1}(c) \subset \R^{n+2}$ if $c>0$, $Q^{n+1}(c) = \R^{n+1}$ if $c=0$ and $Q^{n+1}(c) \subset \R^{n+2}_1$ if $c<0$, we can view $\nu$ as a map from $S$ to $S^{n+1}(1)$ if $c>0$, to $S^n(1)$ if $c=0$ and to $S^{n+1}_1(1)$ if $c<0$. This map is called the \textit{Gauss map} of the hypersurface.
\end{definition}

\subsection{Local Description in Lorentzian Space Forms}

The flat and non-flat cases are treated separately. In both theorems, a distinction is made between the submanifolds with null second fundamental form and the other marginally trapped submanifolds. 

\begin{theorem} \cite{localspace} \label{T:LocalMinkowski}
\begin{itemize}
\item[(1)]Let $\Omega$ be an open domain of $\mathbb{R}^n$ and $\tau\in C^2(\Omega)$, such that $\Delta \tau$ is never zero, where $\Delta$ is the Laplacian of $\mathbb{R}^n$. Then, the immersion $\phi:\Omega\rightarrow\mathbb{R}^{n+2}_1$, defined by
\begin{align*}
\phi(x)=(\tau(x), x, \tau(x)),
\end{align*}
is flat and its second fundamental form is given by
\begin{align*}
h(X,Y)=\mathrm{Hess}_{\tau}(X,Y)(1,0,\ldots, 0, 1).
\end{align*}
In particular, $\phi$ has null second fundamental form and is therefore marginally trapped.

Conversely, any $n$-dimensional submanifold of $\mathbb R^{n+2}_1$ with null second fundamental form is locally congruent to the image of such an immersion. 

\item[(2)]Let $\varphi$ be an immersion of class $C^4$ of an $n$-dimensional manifold $S$ into $\mathbb{R}^{n+1}$ and denote by $\nu$ the Gauss map of $\varphi$. Assume that $\varphi$ admits $p\geq2$ distinct, non-vanishing principal curvatures $\kappa_1,\ldots,\kappa_p$ with multiplicities $m_1,\ldots,m_p$ respectively and denote by $\tau_1,\ldots,\tau_{p-1}$ the $p-1$ roots of the polynomial
\begin{align*}
P(\tau)=\sum_{i=1}^{p}m_i\prod_{j\neq i}^p(\kappa_{j}^{-1}-\tau).
\end{align*}
Then, the $p-1$ immersions $\phi_i:S\rightarrow \mathbb{R}^{n+2}_1$, defined by
\begin{align*}
\phi_i(x) = (\tau_i(x), \, \varphi(x)+\tau_{i}(x)\nu(x))
\end{align*}
for $i \in \{1,\ldots,p-1\}$, are marginally trapped.

Conversely, any $n$-dimensional marginally trapped submanifold of $\mathbb{R}^{n+2}_1$ whose second fundamental form is not null is locally congruent to the image of such an immersion. 
\end{itemize}
\end{theorem}

\begin{remark} Recall that the zero vector is by definition spacelike. The condition that $\Delta\tau$ is nowhere vanishing in case (1) of the above theorem ensures that the mean curvature of the corresponding immersion $\phi$ is nowhere zero. However, the submanifolds described in case (2) do include examples with zero mean curvature. For example, if $\varphi$ is minimal, then
$$ P(0) = \frac{m_1\kappa_1 + \ldots + m_p\kappa_p}{\kappa_1 \ldots \kappa_p} = 0 $$
and the immersion defined by $\phi(x) = (0,\varphi(x))$ is minimal. Similar remarks apply to all the following theorems in this section.
\end{remark}

%\newpage

\begin{theorem} \cite{localspace}. \label{T:localspaceform}
\begin{itemize}
\item[(1)] Let $\Omega$ be an open domain of $S^n(1) \subset \R^{n+1}_1$ (respectively, $H^n(-1)\subset \R^{n+1}_2$) and $\tau\in C^2(\Omega)$, such that $\Delta \tau$ is never zero, where $\Delta$ is the Laplacian of $S^n(1)$ (respectively $H^n(-1)$). Then the immersion $\phi:S^n(1)\rightarrow S^{n+2}_1(1)$ (respectively, $H^n(-1)\rightarrow H^{n+2}_1(-1)$), defined by 
\begin{align*}
\phi(x) = (\tau(x), x, \tau(x)),
\end{align*}
is flat and its second fundamental form is given by
\begin{align*}
h(X,Y) = \mathrm{Hess}_{\tau}(X,Y)(1,0,\ldots, 0, 1).
\end{align*}
In particular, $\phi$ has null second fundamental form and is therefore marginally trapped.

Conversely, any $n$-dimensional spacelike submanifold with null second fundamental form is locally congruent to the image of such an immersion.

\item[(2)] Let $\varphi$ be an immersion of class $C^4$ of an $n$-dimensional manifold $S$ into $S^{n+1}(1)$ (respectively, $H^{n+1}(-1)$) and denote by $\nu$ the $S^{n+1}(1)$-valued (respectively $S^{n+1}_1(1)$-valued) Gauss map of $\varphi$. Assume that $\varphi$ admits $p \geq 2$ distinct, non-vanishing principal curvatures $\kappa_1, \ldots,\kappa_p$ with multiplicities $m_1,\ldots,m_p$ respectively and denote by $\tau_1,\ldots,\tau_{p-1}$ the $p-1$ roots of the polynomial
\begin{align*}
P(\tau)=\sum_{i=1}^{p}m_i\prod_{j\neq i}^p(\kappa_{j}^{-1}-\tau).
\end{align*}
Then, the $p-1$ immersions $\phi_i:S\rightarrow S^{n+2}_1(1)$ (respectively $H^{n+2}_1(-1)$), defined by
\begin{align*}
\phi_i(x) = (\tau_i(x), \, \varphi(x)+\tau_{i}(x)\nu(x))
\end{align*}
for $i \in \{1,\ldots,p-1\}$ are marginally trapped.

Conversely, any $n$-dimensional marginally trapped submanifold of $S^{n+2}_1(1)$ (respectively, of $H^{n+2}_1(-1)$), whose second fundamental form is not null is locally congruent to the image of such an immersion.
\end{itemize}
\end{theorem}

\subsection{Local Description in Robertson-Walker Spacetimes}

Local representation formulas for codimension two marginally trapped submanifolds of the Lorentzian product of the real line and a space form, corresponding to the warping function $f$ of $L^n_1(f,c)$ being constant, were proven in \cite{localspace}. These results were further generalized to marginally trapped submanifolds of spaces of type $L^n_1(f,c)$ for arbitrary positive smooth $f$ in \cite{localRobertson}. We state the results below.

\begin{theorem} \cite{localspace}
\begin{itemize}
\item[(1)]There are no $n$-dimensional submanifolds of $\mathbb{R}_1\times S^{n+1}(1)$ with null second fundamental form.  
\item[(2)]
Let $\varphi$ be an immersion of class $C^4$ of an $n$-dimensional manifold $S$ into $S^{n+1}(1)$. Denote by $\nu$ the Gauss map of $\varphi$ and by $\kappa_1, \ldots, \kappa_p$ its $p$ distinct principal curvatures with multiplicities $m_1,\ldots,m_p$ respectively. Then, the polynomial
\begin{align*}
P(\tau) = \sum_{i=1}^{p}m_i(\kappa_{i}\tau+1)\prod_{j\neq i}^{p}(\tau-\kappa_j)
\end{align*}
has exactly $p-1$ distinct roots $\tau_i$ if $\varphi$ is minimal and $p$ distinct roots otherwise. Moreover, the $p-1$ or $p$ immersions $\phi_i: S \rightarrow \mathbb{R}_1 \times  S^{n+1}(1)$ defined by
\begin{align*}
\phi_i = \left(\cot^{-1}\tau_i,\frac{\tau_{i}\varphi+\nu}{\sqrt{1+\tau_{i}^2}}\right)
\end{align*}
are marginally trapped. 

Conversely, any $n$-dimensional marginally trapped submanifold of $\mathbb{R}_1 \times  S^{n+1}(1)$ is locally congruent to the image of such an immersion. 
\end{itemize}
\end{theorem}

\begin{theorem} \cite{localspace}
\begin{itemize}
\item[(1)]
There are no $n$-dimensional submanifolds of $\mathbb{R}_1\times H^{n+1}(-1)$ with null second fundamental form.  
\item[(2)]
Let $\varphi$ be an immersion of class $C^4$ of an $n$-dimensional manifold $S$ into $H^{n+1}(-1)$. Denote by $\nu$ the Gauss map of $\varphi$ and by $\kappa_1, \ldots, \kappa_p$ its $p$ distinct principal curvatures with multiplicities $m_1,\ldots,m_p$ respectively. Denote by $\tau_1,\ldots,\tau_q$ the different roots of the polynomial 
\begin{align*}
P(\tau) = \sum_{i=1}^{p}m_i(\kappa_{i}\tau-1)\prod_{j\neq i}^{p}(\tau-\kappa_j),
\end{align*}
satisfying $|\tau_i|>1$ for $i \in \{1,\ldots,q\}$. Then the $q \leq p$ immersions $\phi_i: S \rightarrow \mathbb{R}_1\times H^{n+1}(-1)$ defined by 
\begin{align*}
\phi_i = \left(\coth^{-1}\tau_i, \frac{\tau_{i}\varphi+\nu}{\sqrt{\tau_{i}^2-1}}\right)
\end{align*}
are marginally trapped. 

Conversely, any $n$-dimensional marginally trapped submanifold of $\mathbb{R}_1\times H^{n+1}(-1)$ is locally congruent to the image of such an immersion. 
\end{itemize}
\end{theorem}

Now consider the space $L^{n+2}_1(f,c)=I\times_f Q^{n+1}(c)$, where $c\in\left\{-1,0,1 \right\}$ as defined in Section \ref{sec:preliminaries}. It is useful to introduce following notation:
\begin{align*}
(\cos_{c} t, \sin_{c} t)=\begin{cases}
    (\cos t, \sin t) & \text{if $c=1$},\\
    (1, t) & \text{if $c=0$},\\
    (\cosh t, \sinh t) & \text{if $c=-1$}.
  \end{cases}
\end{align*}
The local description of codimension two marginally trapped submanifolds of $L^n_1(f,c)$ is as follows.

\begin{theorem} \cite{localRobertson} \label{T:localRobertson} 
\begin{itemize}
\item[(1)] Let $\phi: S \to L^{n+2}_1(f,c)$ be an immersion of an $n$-dimensional submanifold with null second fundamental form. Then there are two possibilities: 
\begin{itemize}
\item[(i)] the immersion takes the form
$$\phi = (t_0,\varphi),$$ 
where $t_0 \in I$ is a constant and $\varphi$ defines a totally umbilical hypersurface of $Q^{n+1}(c)$ with a mean curvature vector of constant length $|f'(t_0)|$;
\item[(ii)] the function 
$$ \frac{\theta''\cos_c\theta+(\theta')^2c\sin_c\theta}{\theta''\sin_c\theta-(\theta')^2\cos_c\theta}, $$
where $\theta$ is defined in \eqref{eq:theta}, is constant, say $C_0$, and $\phi$ is locally congruent to an immersion of the form 
$$\phi = \left(\tau, \cos_c(\theta\circ\tau)\varphi+\sin_c(\theta\circ\tau)\nu\right),$$ 
where $\tau: S \to \R$ is a real function of class $C^2$  and $\varphi$ defines a totally umbilical hypersurface of $Q^{n+1}(c)$ with Gauss map $\nu$ and with a mean curvature vector of constant length $|C_0|$.
\end{itemize}
\item[(2)]Let $\varphi:S\rightarrow Q^{n+1}(c)$, with $c\in\left\{-1,0,1 \right\}$, be an immersed hypersurface of class $C^4$ with Gauss map $\nu$. Denote by $\kappa_1, \ldots, \kappa_p$ the $p \geq 2$ distinct principal curvatures with multiplicities $m_1,\ldots,m_p$ of $\varphi$. Consider the immersion $\phi:S\rightarrow L^{n+2}_1(f,c)$, defined by
\begin{align*}
\phi = \left(\tau, \cos_c(\theta\circ\tau)\varphi+\sin_c(\theta\circ\tau)\nu\right),
\end{align*} 
where $\tau$ is a real function in $C^2(S)$ and $\theta$ is defined as in equation (\ref{eq:theta}). The immersion $\phi$ is marginally trapped if and only if $\tau: S \to \R$ satisfies
\begin{align*}
n \left(\frac{\de f}{\de t}\circ\tau\right) - \sum_{i=1}^{p}m_i\frac{\kappa_{i}\cos_{c}(\theta\circ\tau)+c\sin_{c}(\theta\circ\tau)}{\cos_{c}(\theta\circ\tau)-\kappa_i\sin_{c}(\theta\circ\tau)}=0.
\end{align*}

Conversely, any marginally trapped codimension two submanifold of $L^{n+2}_1(f,c)$ is locally congruent to the image of such an immersion. 

\end{itemize}
\end{theorem}

%%%%%%%%%%%%%%%%%%%%%%%%%%%%%%%%%%%%%%%%%%%%%%%%%%%%%%%%%%%

\section{Marginally Trapped Surfaces with Positive Relative Nullity}
\label{sec:posrelnul}

Let $S$ be a spacelike submanifold of a pseudo-Riemannian manifold and denote the second fundamental form by $h$. The \textit{relative null space} at a point $p\in S$ is defined by
\begin{align*}
\mathcal{N}_p(S) = \left\{X\in T_p S \ | \ h(X,Y)=0 \text{ for all } Y\in T_p S \right\}.
\end{align*}
The dimension of $\mathcal{N}_p(S)$ is called the \textit{relative nullity of $S$ at $p$} and the submanifold $S$ is said to have \textit{positive relative nullity} if $\dim\mathcal{N}_p(S)>0$ for all $p\in S$.

\subsection{Classification in Lorentzian Space Forms}

The following result classifies all marginally trapped surfaces with positive relative nullity in the Minkowski spacetime $\mathbb{R}^4_1$.

\begin{theorem}\cite{Positiverelativenullity}\label{T:Positive:R}
Up to isometries, there are two families of marginally trapped surfaces with positive relative nullity in the Minkowski spacetime $\mathbb{R}^4_1$:
\begin{itemize}
\item[(1)] a surface parametrized by $\phi(x,y) = \big( f(x),x,y, f(x) \big)$, where  $f$ is an arbitrary differentiable function such that $f''$ vanishes nowhere;
\item[(2)] a surface parametrized by
\begin{multline*} 
\phi(x,y) = \Bigg( \int_0^x r(s)q'(s) \, \de s + q(x)y, \, y\cos x - \int_0^x r(s)\sin s \, \de s, \\ 
y\sin x + \int_0^x r(s)\cos s \, \de s, \, \int_0^x r(s)q'(s) \, \de s + q(x)y \Bigg), 
\end{multline*}
where $q$ and $r$ are  defined on an open interval $I$ containing $0$, such that $q''+q$ vanishes nowhere on $I$. 
\end{itemize}

Conversely,  every  marginally trapped surface with positive relative nullity in the Minkowski spacetime $\mathbb{R}^4_1$ is congruent to an open part of a surface in one of the two families.
\end{theorem}

The next corollary follows immediately from the two explicit parametrizations in Theorem \ref{T:Positive:R} since the first and fourth components are equal to each other in both paramatrizations. However, one can also prove it directly by using the Erbacher-Magid reduction theorem from \cite{Magid}, see for example \cite{Walleghem}.

\begin{corollary}\label{T:NullhyperE}
Every marginally trapped surface with positive relative nullity in $\mathbb{R}^4_1$ is contained in a null hyperplane of $\mathbb{R}^4_1$.
\end{corollary}

The following two theorems give classifications of marginally trapped surfaces with positive relative nullity in the de Sitter and anti-de Sitter spacetimes. 

\begin{theorem}\cite{Positiverelativenullity} \label{T:Positive:S} 
Up to isometries, there are two families of  marginally trapped surfaces with positive relative nullity  in the de Sitter spacetime $S^4_1(1)$:
\begin{itemize}
\item[(1)] a surface parametrized by
\begin{align*} 
\phi(x,y) = \Big( f(x) \cos y, \, \sin x\cos y, \, \sin y, \, \cos x\cos y, \, f(x)\cos y \Big), 
\end{align*}
where $f$ is an arbitrary differentiable function such that $f''+f$ vanishes nowhere;
\item[(2)] a surface parametrized by
\begin{multline*}
\phi(x,y)= \big( p(x), \, \eta_1(x), \, \eta_2(x), \, \eta_3(x), \, p(x) \big) \cos y \\ 
- \Big( b-\int_0^x r(s)p'(s)\de s, \, \xi_1(s), \, \xi_2(s), \, \xi_3(s), \, b-\int_0^x r(s)p'(s)\de s \Big) \sin y,
\end{multline*}
where $b$ is a real number, $p$ and $r$ are defined on an open interval $I$ containing $0$ such that $r$ is non-constant, $\eta=(\eta_1,\eta_2,\eta_3)$ is a unit speed curve in $S^2(1)\subset \mathbb{R}^3$  with geodesic curvature $\kappa_g=r$ and $\xi=(\xi_1,\xi_2,\xi_3)$ is the unit normal of $\eta$ in $S^2(1)$.
\end{itemize}

Conversely, every  marginally trapped surface with positive relative nullity in the de Sitter spacetime $S^4_1(1)$ is congruent to an open part of a surface in one of the two families.
\end{theorem}

\begin{theorem} \label{T:Positive:H}\cite{Positiverelativenullity} 
Up to isometries, there are five families of  marginally trapped surfaces with positive relative nullity  in the anti-de Sitter spacetime $H^4_1(-1)$:
\begin{itemize}
\item[(1)] a surface parametrized by
\begin{align*}
\phi(x,y) = \Big( f(x)\cosh y, \, \cosh x\cosh y, \, \sinh y, \, \sinh x\cosh y, \, f(x)\cosh y \Big),
\end{align*}
where $f$ is a differentiable function such that $f''-f$ vanishes nowhere;
\item[(2)] a surface parametrized by
\begin{align*}
\phi(x,y) = \Big( f(x)\sinh y, \, \cosh y, \, \cos x \sinh y, \, \sin x \sinh y,  \, f(x)\sinh y \Big),
\end{align*}
where $f$ is a differentiable function such that $f''-f$ vanishes nowhere;
\item[(3)] a surface parametrized by
\begin{align*}
\phi(x,y) = \( x^2 e^y, \, \frac{3e^y}{2}-2\sinh y, \, e^y-2\sinh y, \, x e^y, \, x^2e^y -\frac{e^y}{2} \);
\end{align*}
\item[(4)] A  surface parametrized by
\begin{align*} 
\phi(x,y) = \( \sinh y-\frac{x^2e^y}{2}-e^y, \, f(x)e^y, \, xe^y, \, f(x)e^y, \, \sinh y-\frac{x^2e^y}{2} \),
\end{align*}
where $f$ is a differentiable function such that $f''$ vanishes nowhere;
\item[(5)] A surface parametrized by
\begin{multline*}
\phi(x,y) = \big( p(x), \, \eta_1(x), \, \eta_2(x), \, \eta_3(x), \, p(x) \big) \cosh y \\ 
-\( b-\int_0^x r(s)p'(s) \, \de s, \, \xi_1(s), \, \xi_2(s), \, \xi_3(s), \, b-\int_0^x r(s)p'(s) \, \de s \)\sinh y,
\end{multline*} 
where $b$ is a real number, $p$ and $r$ are defined on an open interval $I$ containing $0$ such that $r$ is non-constant, $\eta=(\eta_1,\eta_2,\eta_3)$ is a unit speed curve in $H^2(-1) \subset \mathbb{R}_1^3$  with geodesic curvature $\kappa_g=r$ and $\xi=(\xi_1,\xi_2,\xi_3)$ is a unit normal of $\eta$ in $H^2(-1)$. 
\end{itemize}

Conversely, every  marginally trapped surface with positive relative nullity in the anti-de Sitter spacetime $H^4_1(-1)$ is congruent to an open part of a surface in one of the five families.
\end{theorem}

Also for de Sitter and anti-de Sitter spacetimes, we have corollaries similar to Corollary \ref{T:NullhyperE}.
 
\begin{corollary}\label{T:NullhyperS}
Every marginally trapped surface with positive relative nullity in $S^4_1(1)$, respectively $H^4_1(-1)$, is contained in a null hyperplane of $\R^5_1$, respectively $\R^5_2$.
\end{corollary}

\subsection{Classification in Robertson-Walker Spacetimes}

It turns out that marginally trapped surfaces with positive relative nullity in Robertson-Walker spacetimes (of non-constant sectional curvature) do not exist.

\begin{theorem}\cite{Robertson}
Let $L^4_1(f,c)$ be a Robertson-Walker spacetime which contains no open subsets of constant sectional curvature. Then $L^4_1(f,c)$ does not admit any marginally trapped surfaces with positive relative nullity. 
\end{theorem}

Note that an open subset of $L^4_1(f,c)$ of constant sectional is isometric to an open part of a Lorentzian space form, so marginally trapped surfaces with positive relative nullity in such a subset are classified in Theorem~\ref{T:Positive:R}, Theorem~\ref{T:Positive:S} and Theorem~\ref{T:Positive:H}. 

%%%%%%%%%%%%%%%%%%%%%%%%%%%%%%%%%%%%%%%%%%%%%%%%%%%%%%%%%%%

\section{Marginally Trapped Surfaces with Parallel Mean Curvature Vector Field}
\label{sec:pmc}

A submanifold $S$ of a pseudo-Riemannian manifold is said to have \textit{parallel mean curvature vector field} if $\nabla_X^{\perp}H=0$ for every vector field $X$ tangent to $S$ and to have \textit{parallel second fundamental form} or to be \textit{parallel} for short, if $\overline{\nabla}_X h=0$ for every vector field $X$ tangent to $S$. Here, $\overline{\nabla}$ is the connection of Van der Waerden-Bortolotti, defined by
$$ (\overline{\nabla}_X h)(Y,Z) = \nabla^{\perp}_X h(Y,Z) - h(\nabla_X Y,Z) - h(Y,\nabla_X Z) $$
for all vector fields $X$, $Y$ and $Z$ tangent to $S$. It is easy to see that a parallel submanifold has parallel mean curvature vector field, but the converse is not necessarily true.

In this section we give the classification of marginally trapped surfaces with parallel mean curvature vector field in four-dimensional Lorentzian space forms, which was proven in \cite{Parallelmeancurvature}. These surfaces are then related to marginally trapped surfaces with a 1-type Gauss map, as shown in \cite{Gaussmap}. The classification of marginally trapped surfaces with parallel mean curvature vector field uses the notion of a \textit{light cone}. 

\begin{definition}
The light cone $\mathcal{LC}$ in Minkowski spacetime $\mathbb{R}^4_1$ is
\begin{align*}
\mathcal{LC} = \left\{ x \in \mathbb{R}^4_1 \ | \ \langle x,x \rangle = 0 \right\}. 
\end{align*}
We can see $\mathcal{LC}$ as a submanifold of de Sitter spacetime, respectively anti-de Sitter spacetime, by using the following natural embeddings:
\begin{align*}
& \mathcal{LC} \to S^4_1(1) \subset \mathbb{R}^5_1: x \mapsto (x,1), \\
& \mathcal{LC} \to H^4_1(-1) \subset \mathbb{R}^5_2: x \mapsto (1,x).
\end{align*}
\end{definition}

The following propositions show that marginally trapped surfaces with parallel mean curvature vector field arise naturally in the light cones of four-dimensional Lorentzian space forms.

\begin{proposition}\cite{Parallelmeancurvature}\label{P:Lightcones}
Let $S$ be a marginally trapped surface in $Q^4_1(c)$, with $c\in\left\{-1,0,1 \right\}$. If S lies in $\mathcal{LC}\subset Q^4_1(c)$, then $S$ has constant Gaussian curvature $c$ and has parallel mean curvature vector field in $Q^4_1(c)$.
\end{proposition}

\begin{proposition}\cite{Parallelmeancurvature}\label{P:ExistenceLightcones}
Let $\lambda$ be a positive solution of the differential equation
\begin{equation} \label{eq:laplacian_lambda}
\lambda(\lambda_{xx}+\lambda_{yy})-\lambda_x^2-\lambda_y^2-2c\lambda = 0
\end{equation}
for $c \in \{-1,0,1\}$ on a simply connected domain $U\subset\mathbb{R}^2$. Then $S=(U, \lambda^{-1}(\de x^2+\de y^2))$ is a surface of constant curvature $c$. Moreover, there exists a marginally trapped isometric immersion $\phi:S\rightarrow Q^4_1(c)$ with parallel mean curvature vector field such that $\phi$ lies in the light cone $\mathcal{LC}\subset Q^4_1(1)$.
\end{proposition}

\begin{remark}
Note that equation \eqref{eq:laplacian_lambda} can be rewritten as $\Delta\lambda = 2c\lambda$, where $\Delta = \lambda(\partial_x^2+\partial_y^2)-\lambda_x\partial_x-\lambda_y\partial_y$ is the Laplacian of $S=(U, \lambda^{-1}(\de x^2+\de y^2))$.
\end{remark}

%%%%%%%%%%%%%%%%%%%%%%%%%%%%%%%%%%%%%%%%%%%%%%%%%%%%%%%%%%%

\subsection{Classification in Lorentzian Space Forms}

The following theorem classifies all marginally trapped surfaces with parallel mean curvature in $\mathbb{R}^4_1$, $S^4_1(1)$ and $H^4_1(-1)$.

\begin{theorem}\cite{Parallelmeancurvature}\label{T:ParallelMeanE}
Let $S$ be a marginally trapped surface with parallel mean curvature vector field in the Minkowski spacetime $\mathbb{R}^4_1$. Then $S$ is one congruent to of the following six types of surfaces:
\vskip.05in
\begin{itemize}
\item[(1)]  a flat parallel surface given by 
$$ \phi(x,y)=\frac{1}{2}\Big((1-b)x^2+(1+b)y^2,(1-b)x^2+(1+b)y^2,2x,2y\Big) $$
for some $b \in \R$;
\vskip.05in
\item[(2)]  a flat parallel surface given by $\phi(x,y)=a\big(\cosh x, \sinh x,\cos y,\sin y\big)$, with $a>0$;
\vskip.05in
\item[(3)]a flat surface given by 
\begin{align*}
\phi(x,y)=(f(x,y), x, y, f(x, y)),
\end{align*}
where $f$ is a smooth function on $S$ such that $\Delta f=a$ for some nonzero real number $a$;
\vskip.05in
\item[(4)] a non-parallel flat surface lying in the light cone $\mathcal{LC}$;
\vskip.05in
\item[(5)]  a non-parallel surface lying in the de Sitter spacetime $S^3_1(c)$ for some $c>0$ such that the  mean curvature vector field $H'$ of $S$ in $S^3_1(c)$ satisfies $\<H',H'\>=-c$;
\vskip.05in
\item[(6)]  a non-parallel surface lying in the hyperbolic space $H^3(c)$ for some $c<0$ such that the mean curvature vector field $H'$ of $S$ in $H^3(c)$ satisfies $\<H',H'\>=-c$.
\end{itemize}
Conversely, every surface of types $(1)$--$(6)$  above gives rise to a marginally trapped surface with parallel mean curvature vector in $\mathbb{R}^4_1$. 
\end{theorem}

\begin{theorem}\cite{Parallelmeancurvature}\label{T:Parallel:S} Let  $S$ be a marginally trapped  surface  with parallel mean curvature vector field in the de Sitter spacetime $S^4_1(1)$.  Then $S$ is congruent to one of the following eight types of surfaces:
\vskip.05in
\begin{itemize}
\item[(1)] a parallel surface of Gaussian curvature $1$ given by
\begin{align*}
\phi(x,y)=\Big(1,\sin x,\cos x \cos y,\cos x\sin y,1\Big)\mbox{;}
\end{align*}
\item[(2)] a flat parallel surface defined by 
\begin{align*}
\phi(x,y)=\frac{1}{2}\big(2x^2-1,2x^2-2,2x,\sin 2y,\cos 2y\big)\mbox{;}
\end{align*}
\item[(3)] a flat parallel surface defined by 
\begin{align*}
\phi(x,y)=\Bigg(\frac{b}{\sqrt{4-b^2}},\frac{\cos (\sqrt{2-b} \, x)}{\sqrt{2-b}},\frac{\sin (\sqrt{2-b} \, x)}{\sqrt{2-b}},\frac{\cos (\sqrt{2+b} \, y)}{\sqrt{2+b}},\frac{\sin (\sqrt{2+b}\, y)}{\sqrt{2+b}} \Bigg),
\end{align*}
with $|b|<2$;
\vskip.05in
\item[(4)] a flat parallel surface defined by 
$$
\phi(x,y)=\Bigg(\frac{\cosh (\sqrt{b-2} \, x)}{\sqrt{b-2}},\frac{\sinh (\sqrt{b-2} \, x)}{\sqrt{b-2}},\frac{\cos (\sqrt{2+b} \, y)}{\sqrt{2+b}},\frac{\sin (\sqrt{2+b} \, y)}{\sqrt{2+b}}, \frac{b}{\sqrt{b^2-4}} \Bigg),
$$
with $b>2$;
\vskip.05in
\item[(5)] a surface of constant curvature one given by 
\begin{align*}
\phi(x,y)=(f(x,y), \cos x, \sin x\cos y, \sin x\sin y, f(x, y)),
\end{align*}
where $f$ is a smooth function satisfying $\Delta f=a$ for some nonzero real number $a$;
\vskip.05in
\item[(6)] a non-parallel surface of curvature one, lying in the light cone $\mathcal {LC}\subset S^4_1(1)$;
\vskip.05in
\item[(7)] a non-parallel surface lying in $S^4_1(1)\cap S^4_1(x_0,c)$, with $x_0\ne 0$ and $c>0$, such that the mean curvature vector field $H'$ of $S$ in $S^4_1(1)\cap S^4_1(x_0,c)$ satisfies $\<H',H'\>=1-c$.
\vskip.05in
\item[(8)] a non-parallel surface lying in $S^4_1(1)\cap H^4(x_0,c)$, with $x_0\ne 0$ and $c<0$, such that the mean curvature vector field $H'$ of $S$ in $S^4_1(1)\cap H^4(x_0,c)$ satisfies $\<H',H'\>=1-c$.
\end{itemize}
Conversely, every surface of types $(1)$--$(8)$ above gives rise to a marginally trapped surface with parallel mean curvature vector in $S^4_1(1)$. 
\end{theorem}

\begin{theorem}\cite{Parallelmeancurvature}\label{T:Parallel:H} Let  $S$ be a marginally trapped  surface  with parallel mean curvature vector field in the anti-de Sitter spacetime $H^4_1(-1)$.  Then, $S$ is congruent to one of the following eight types of surfaces:
\begin{itemize}
\item[(1)]  a parallel surface of Gaussian curvature $-1$ given by
\begin{align*}
\phi(x,y)=\big(1,\cosh x \cosh y,\sinh x,\cosh x\sinh y,  1\big) \mbox{;}
\end{align*}
\vskip.05in
\item[(2)]   a flat parallel surface defined by 
\begin{align*}
\phi(x,y)=\frac{1}{2}\big(2x^2+2,\cosh 2y, 2x,\sinh 2y,2x^2+1\big) \mbox{;}
\end{align*}
\item[(3)]  a flat parallel surface defined by 
\begin{align*}  
&\phi(x,y) = \Bigg(\frac{\cosh (\sqrt{2\hskip-.02in -\hskip-.02in b} \, x)}{\sqrt{2-b}},\frac{\cosh (\sqrt{2\hskip-.02in +\hskip-.02in b} \, y)}{\sqrt{2+b}}, \frac{\sinh (\sqrt{2\hskip-.02in -\hskip-.02in b} \, x)}{\sqrt{2-b}},\frac{\sinh (\sqrt{2\hskip-.02in +\hskip-.02in b} \, y)}{\sqrt{2+b}},\frac{b}{\sqrt{4\hskip-.02in -\hskip-.02in b^2}} \Bigg), 
\end{align*}
with $|b|<2$;
\vskip.05in
\item[(4)] a flat parallel surface defined by 
\begin{align*}  
&\phi(x,y) = \Bigg(\frac{b}{\sqrt{b^2-4}},\frac{\cosh (\sqrt{b+2} \, y)}{\sqrt{b+2}},\frac{\sinh (\sqrt{b+2} \, y)}{\sqrt{b+2}}, \frac{\cos (\sqrt{b-2} \, x)}{\sqrt{b-2}},\frac{\sin (\sqrt{b-2} \, x)}{\sqrt{b-2}} \Bigg),
\end{align*} 
with $b>2$;
\vskip.05in
\item[(5)] a surface of constant curvature $-1$ given by
\begin{align*}
\phi(x,y)=(f(x,y), \cosh x, \sinh x\cos y, \sinh x\sin y, f(x,y)),
\end{align*}
where $f$ is a smooth function satisfying $\Delta f=a$ for some nonzero real number $a$;
\vskip.05in
\item[(6)] a non-parallel surface with curvature $-1$, lying in the light cone $\mathcal {LC}\subset H^4_1(-1)$;
\vskip.05in
\item[(7)] a non-parallel surface lying in $H^4_1(-1)\cap S^4_2(x_0,c)$, with $x_0\ne 0$ and $c>0$, such that the mean curvature vector field $H'$  in $H^4_1(-1)\cap S^4_2(x_0,c)$ satisfies $\<H',H'\>=-1-c$;
\vskip.05in
\item[(8)]  a non-parallel surface lying  in $H^4_1(-1)\cap H_1^4(x_0,c)$, with $x_0\ne 0$ and $c<0$, such that the mean curvature vector field $H'$ in $H^4_1(-1)\cap H_1^4(x_0,c)$ satisfies $\<H',H'\>=-1-c$.
\end{itemize}

Conversely, every surface of types $(1)$--$(8)$  above gives rise to a marginally trapped surface with parallel mean curvature vector in $H^4_1(-1)$. 
\end{theorem}

\subsection{Finite Type Gauss Map}

In \cite{Gaussmap}, marginally trapped surfaces with parallel mean curvature vector field were related to marginally trapped surfaces with a 1-type Gauss map. This notion of Gauss map is a little different from the one we used in Section \ref{sec:local_description}, so we start by recalling the definition.

\begin{definition}
Let $S \to \R^4_1$, respectively $S \to S^4_1(1) \subset \R^5_1$ or $S \to H^4_1(-1) \subset \R^5_2$, be a spacelike oriented surface in a four-dimensional Lorentzian space form. For any $p \in S$, we denote $\nu(p) = e_1 \wedge e_2$, where $(e_1,e_2)$ is a positively oriented orthonormal basis of $T_pS$. Then $\nu$ is a map from $S$ to $\bigwedge^2 \mathbb R^4_1 \cong \mathbb R^6_3$, respectively $\bigwedge^2 \mathbb R^5_1 \cong \mathbb R^{10}_4$ or $\bigwedge^2 \mathbb R^5_2 \cong \mathbb R^{10}_6$, which we call the \textit{Gauss map} of the surface.
\end{definition}

It now makes sense to look at $\Delta \nu$, where $\Delta$ is the Laplacian of $S$ acting on every component of $\nu$ and we have the following definition, which extends the notion of having harmonic Gauss map.

\begin{definition}
An oriented surface $S$ in $\R^4_1$, respectively $S^4_1(1) \subset \R^5_1$ or $H^4_1(-1) \subset \R^5_2$, has \textit{pointwise $1$-type Gauss map of the first kind} if its Gauss map $\nu$ satisfies
$$ \Delta \nu = f \nu $$
for some function $f: S \to \R$.
\end{definition}

The following result is a combination of Proposition 1, Theorem 5 and Proposition 6 from \cite{Gaussmap}.

\begin{theorem} \cite{Gaussmap}
A marginally trapped surface in $\R^4_1$, $S^4_1(1)$ or $H^4_1(-1)$ has pointwise $1$-type Gauss map of the first kind if and only if it has parallel mean curvature vector field. In particular, a marginally trapped surface in $\R^4_1$ has harmonic Gauss map if and only if it has parallel mean curvature vector field and it is flat, while a marginally trapped surface in $S^4_1(1)$ or $H^4_1(-1)$ cannot have harmonic Gauss map.
\end{theorem}

%%%%%%%%%%%%%%%%%%%%%%%%%%%%%%%%%%%%%%%%%%%%%%%%%%%%%%%%%%%

\section{Marginally Trapped Surfaces in $\mathbb{R}^4_1$ Which Are Invariant under a $1$-Parameter Group of Isometries}\label{S:Invariant}
\label{sec:invariant}

Marginally trapped surfaces in the Minkowski spacetime $\mathbb{R}^4_1$, satisfying the additional condition of being invariant under the action of a $1$-parameter subgroup $G$ of the isometry group of $\R^4_1$, are studied in \cite{Boost,Rotation,Screw}. The main results are the classifications of boost, rotation and screw invariant marginally trapped surfaces in $\R^4_1$ respectively, which we discuss in this section.
 
\subsection{Boost Invariant Marginally Trapped Surfaces}

One says that a spacelike surface $S$ in $\mathbb{R}^4_1$ is invariant under boosts if it is invariant under the following group of linear isometries of $\mathbb{R}^4_1$:
\begin{align*}
\left \{ \left. B_\theta = \begin{pmatrix}
 \cosh\theta & \sinh\theta & 0 & 0 \\ 
 \sinh\theta & \cosh\theta & 0 & 0 \\ 
 0 & 0 & 1 & 0 \\ 
 0 & 0 & 0 & 1 
\end{pmatrix}\ \right| \ \theta\in\mathbb{R}  \right \}.
\end{align*}
This means that $B_\theta S=S$ for all $\theta\in\mathbb{R}$. A boost invariant surface $S$ has an open and dense subset $\Sigma_\alpha$ which can be parametrized by
$$ \phi: I \times \R \to \R^4_1: (s,\theta) \mapsto B_\theta \begin{pmatrix} \alpha_1(s) \\ 0 \\ \alpha_3(s) \\ \alpha_4(s) \end{pmatrix}, $$
where $\alpha: I \subset \mathbb{R}\rightarrow \left\{(x_1, x_2, x_3, x_4)\in \mathbb{R}^4_1 \ | \ x_1>0, \ x_2=0 \right\}$ is a spacelike curve parametrized by arc length. The next theorem classifies all boost invariant marginally trapped surfaces in $\mathbb{R}^4_1$.

\begin{theorem}\cite{Boost}\label{T:Boost}
Let $S$ be a boost invariant marginally trapped surface in the Minkowski spacetime $\mathbb{R}^4_1$. Then, $S$ is locally congruent to a surface $\Sigma_{\alpha}$ whose profile curve $\alpha(s)=(\alpha_1(s), 0, \alpha_3(s), \alpha_4(s))$ is described as follows. Take a positive smooth function $\alpha_1: I \subset \R \to \R$ such that the function $\rho$, defined by 
\begin{align*}
\rho(u)=\frac{1+({\alpha_1}'(u))^2+\alpha_1(u){\alpha_1}''(u)}{\alpha_1(u)},
\end{align*}
is never zero.
Choose a function $\epsilon:I\rightarrow\left\{-1,1\right\}$ such that $\epsilon\rho$ is smooth and define the functions $\alpha_3$ and $\alpha_4$ by
\begin{align*}
\alpha_3(s) = \int_{s_0}^s \sqrt{1+({\alpha_1}'(t))^2} \cos\xi(t) \, \de t, \quad
\alpha_4(s) = \int_{s_0}^s \sqrt{1+({\alpha_1}'(t))^2}\sin\xi(t) \, \de t
\end{align*}
for some $s_0\in I$, where $ \displaystyle{\xi(t) = \int_{s_0}^t \frac{\epsilon(u)\rho(u)}{1+({\alpha_1}'(u))^2} \, \de u} $.
\end{theorem}

\subsection{Rotation Invariant Marginally Trapped Surfaces}

One says that a spacelike surface $S$ in $\mathbb{R}^4_1$ is invariant under (spacelike) rotations if it is invariant under the following group of linear isometries of $\mathbb{R}^4_1$:
\begin{align*}
\left \{ \left. R_\theta=\begin{pmatrix}
 1 & 0  &0  &0 \\ 
 0 & 1  &0  &0 \\ 
 0 & 0 & \cos\theta & -\sin\theta \\ 
 0 & 0 & \sin\theta & \cos\theta
\end{pmatrix}\ \right| \ \theta\in\mathbb{R} \right\},
\end{align*}
which means that $R_\theta S=S$ for all $\theta\in\mathbb{R}$. A rotation invariant surface $S$ has an open and dense subset $\Sigma_{\alpha}$ which can be parametrized by
$$ \phi: I \times \R \to \R^4_1: (s,\theta) \mapsto R_\theta \begin{pmatrix} \alpha_1(s) \\ \alpha_2(s) \\ \alpha_3(s) \\ 0 \end{pmatrix}, $$
where $\alpha: I\subset\R \to \{(x_1, x_2, x_3, x_4) \in \mathbb{R}^4_1\ | \ x_3>0, \ x_4=0 \}$ is a spacelike curve parametrized by arc length. Rotation invariant marginally trapped surfaces in $\mathbb{R}^4_1$ are completely classified in the following theorem.

\begin{theorem}\cite{Rotation}\label{T:Rotation}
Let $S$ be a rotation invariant marginally trapped surface in Minkowski spacetime $\mathbb{R}^4_1$. Then, $S$ is locally congruent to a surface $\Sigma_{\alpha}$ whose profile curve $\alpha(s)=(\alpha_1(s),\alpha_2(s),\alpha_3(s),0)$ is described in one the following two cases.
\begin{itemize}
\item[(1)] Given a smooth function $\tau: I \subset (0, +\infty) \rightarrow \mathbb{R}$, such that $\tau(s)+s\tau'(s)$ is never zero,  choose a function $\epsilon:I\rightarrow\left\{-1,1\right\}$ such that $\epsilon\tau$ is also smooth. Define the functions $\alpha_1,\alpha_2,\alpha_3$ by
\begin{align*}
\alpha_1(s)=\int_{s_0}^s \epsilon(t)\tau(t) \, \de t, \quad \alpha_2(s)=\int_{s_0}^s \tau(t) \, \de t, \quad \alpha_3(s)=s
\end{align*}
for some $s_0\in I$.  
\item[(2)] Given a smooth positive function $\alpha_3: I\subset\R\to\R$, such that $1-(\alpha_3')^2-\alpha_3\alpha_3''$ is never zero, define the functions $\alpha_1$ and $\alpha_2$ by
\begin{align*}
& \alpha_1(s) = \int_{s_0}^s\left( \sinh\xi(t)-\alpha'_3(t)\cosh\xi(t) \right)\de t,\\
& \alpha_2(t) = \int_{s_0}^s \left( \cosh\xi(t)-\alpha'_3(t)\sinh\xi(t) \right) \de t
\end{align*}
for some $s_0\in I$, where $ \displaystyle{\xi(t) = \int_{s_0}^t \frac{\de u}{\alpha_3(u)}} $.
\end{itemize}
\end{theorem}

\subsection{Screw Invariant Marginally Trapped Surfaces}

One says that a spacelike surface $S$ in $\mathbb{R}^4_1$ is screw invariant if it is invariant under the following group of linear isometries of $\mathbb{R}^4_1$:

\begin{align*}
\left\{ \left. S_\theta = \begin{pmatrix}
1 & \theta^2 & \sqrt{2}\,\theta & 0 \\ 
0 & 1 & 0 & 0 \\ 
0 & \sqrt{2}\,\theta & 1 & 0 \\ 
0 & 0 & 0 & 1 
\end{pmatrix} \ \right| \ \theta\in\mathbb{R} \right\},
\end{align*}
where the matrices are written with respect to the ordered basis $(k, l, e_3 ,e_4)$, with $k=(1, 1, 0, 0)/\sqrt{2}$, $l=(1, -1, 0, 0)/\sqrt{2}$, $e_3=(0,0,1,0)$ and $e_4=(0,0,0,1)$.  Note that $k$ and $l$ are null vectors with $\langle k,l \rangle = -1$. A screw invariant surface $S$ is contained in $\mathcal{R}^+\cup\mathcal{R}^{-}$, defined as
\begin{align*}
&\mathcal{R}^+ = \left\{ x_kk+x_ll+x_3e_3+x_4e_4 \in \mathbb{R}^4_1 \ | \ x_k>0\right\}, \\
&\mathcal{R}^- = \left\{ x_kk+x_ll+x_3e_3+x_4e_4 \in \mathbb{R}^4_1 \ | \ x_k<0\right\}.
\end{align*}
We will suppose that $S$ is contained in $\mathcal{R}^+$, as similiar results can be obtained for surfaces lying in $\mathcal{R}^-$.  A screw invariant surface $S$ in $\mathcal R^+$ has an open and dense subset $\Sigma_{\alpha}$ which can be parametrized by
$$ \phi: I \times \R \to \R^4_1: (s,\theta) \mapsto S_\theta \begin{pmatrix} \alpha_k(s) \\ \alpha_l(s) \\ 0 \\ \alpha_4(s) \end{pmatrix}, $$
where $\alpha: I\subset\R \to \{ x_kk+x_ll+x_3e_3+x_4e_4 \in \mathbb{R}^4_1\ | \ x_k>0, \ x_3=0 \}$ is a spacelike curve parametrized by arc length. Note that $\alpha_k$, $\alpha_l$, $0$ and $\alpha_4$ are the coordinate functions of $\alpha$ with respect to the basis $(k, l, e_3 ,e_4)$. The next theorem classifies all screw invariant marginally trapped surfaces in $\mathcal{R}^+$. Similar results can also be obtained for surfaces lying in $\mathcal{R}^-$.

\begin{theorem}\cite{Screw}\label{T:Screw}
Let $S$ be a screw invariant marginally trapped surface in $\mathcal{R}^{+}$. Then, $S$ is locally congruent to a surface $\Sigma_\alpha$ whose profile curve $\alpha(s)=\alpha_k(s)k+\alpha_l(s)l+\alpha_4(s)e_4$ is described in one of the following two cases. 
\begin{itemize}
\item[(1)] $\alpha_k$ is a positive constant, $\alpha_4(s) = s + s_0$ for some $s_0 \in \R$ and $\alpha_l$ is a smooth function such that $1+\alpha_k\alpha''_l$ is never zero. 
\item[(2)]Given two functions $\rho:I\subset\mathbb{R}\rightarrow\mathbb{R}$ and $\epsilon:I\rightarrow\left\{-1,1\right\}$, such that $\rho$ and $\epsilon\rho$ are smooth and $\rho$ is never zero, define the functions $\alpha_k, \alpha_l, \alpha_4 :I\rightarrow\mathbb{R}$ as follows:
\begin{align*}
& \alpha_k(s) = \sqrt{\alpha_{k0}+\int_{s_0}^s\exp \xi(t) \, \de t},\\
& \alpha_l(s) = \int_{s_0}^s\frac{\alpha_k(t)}{4\exp \xi(t)}\left(\alpha_{k0}+\int_{s_0}^t\frac{2\alpha_k(u)\rho(u)}{\exp \xi(u)} \, \de u \right)^2\de t - \int_{s_0}^s\frac{2\alpha_k(t)}{\exp \xi(t)} \, \de t, \\
& \alpha_4(s) = \alpha_{40}\alpha_k(s)+\int_{s_0}^s\frac{\exp \xi(t)}{\alpha_k(t)}\left(\int_{s_0}^t\frac{\alpha_k(u)\rho(u)}{\exp \xi(u)} \, \de u \right) \de t,
\end{align*}
with $\alpha_{k0}, \alpha_{40} \in\mathbb{R}$ and $s_0\in I$, where $ \displaystyle{\xi(t) = \int_{s_0}^t\epsilon(u)\rho(u) \, \de u+\xi_0} $ for some $\xi_0 \in \R$.
\end{itemize}
\end{theorem}

%%%%%%%%%%%%%%%%%%%%%%%%%%%%%%%%%%%%%%%%%%%%%%%%%%%%%%%%%%%

\section{Isotropic Marginally Trapped Surfaces}
\label{sec:isotropic}

A complete classification of complete isotropic marginally trapped surfaces in Lorentzian space forms was obtained in \cite{Isotropy}.

\begin{definition} \label{def:isotropic}
An isometric immersion of a Riemannian manifold $S$ into a (pseudo-)Riemannian manifold is called \textit{isotropic} if $\langle h(u,u),h(u,u)\rangle = \lambda(p)$ does not depend on the choice of the unit vector $u\in T_p S$. The~function $\lambda: S\to\R$ is then called the \textit{isotropy function} of the immersion.
\end{definition}

\begin{definition} \label{def:pseudo-umbilical}
An isometric immersion of a Riemannian manifold $S$ into a (pseudo-)Riemannian manifold is called \textit{pseudo-umbilical} if there exists a function $\rho: S \to \R$ such that
$\langle h(X,Y),H \rangle = \rho \langle X,Y \rangle$ for all vector fields $X$ and $Y$ tangent to $S$. 
\end{definition}

Remark that the function $\rho$ in Definition \ref{def:pseudo-umbilical} has to equal $\rho = \langle H,H \rangle$. When the mean curvature vector of a spacelike surface in a four-dimensional Lorenztian manifold is null, the notions of isotropy and pseudo-umbilicity are equivalent. More precisely, the following proposition was proven in \cite{Isotropy}. 

\begin{proposition}\cite{Isotropy} \label{P:Isotropy}
Let $\phi:S\rightarrow M^4_1$ be a marginally trapped surface in a four-dimensional Lorentzian manifold. Then the following assertions are equivalent: 
\begin{itemize}
\item[(1)] $\phi$ is pseudo-umbilical,
\item[(2)] $\phi$ is isotropic,
\item[(3)] $\phi$ is isotropic with isotropy function $0$,
\item[(4)] $\phi$ has null second fundamental form.
\end{itemize}
\end{proposition}
\subsection{Classification in Lorentzian Space Forms}
The following theorem, which unifies Theorem 5.6, Theorem 5.10 and Theorem 5.13 from \cite{Isotropy}, classifies the isotropic marginally trapped surfaces in $\mathbb{R}^4_1, S^4_1(1)$ and $H^4_1(-1)$. Proposition \ref{P:Isotropy} then allows the condition of being isotropic to be replaced by the stronger sounding condition of being isotropic with isotropy function $0$ or by the condition of being pseudo-umbilical. 

\begin{theorem}\cite{Isotropy} \label{T:isotropy}
Let $S$ be a complete connected spacelike surface in $Q^4_1(c)$, with $c \in \{-1,0,1\}$. Then $S$ is an isotropic marginally trapped surface if and only if $S$ has constant Gaussian curvature $c$ and the immersion is congruent to 
$$ \phi : Q^2(c) \to Q^4_1(c) : x \mapsto (\tau(x),x,\tau(x)) $$
where $\tau: Q^2(c) \to \R$ is a smooth function such that $\Delta \tau$ is nowhere zero, where $\Delta$ is the Laplacian of $Q^2(c)$.
\end{theorem}

Remark that, since isotropic marginally trapped surfaces in four-dimensional Lorentzian space forms have null second fundamental form, Theorem \ref{T:isotropy} follows from the first cases of Theorem \ref{T:LocalMinkowski} and Theorem \ref{T:localspaceform}, which were historically proven after Theorem \ref{T:isotropy}.

\subsection{Classification in Robertson-Walker spacetimes}
Isotropic marginally trapped submanifolds in Robertson-Walker spacetimes can now be classified as a corollary of Proposition \ref{P:Isotropy} and the first cases of Theorem \ref{T:localRobertson}.

\begin{theorem}
Let $\phi: S \to L^{4}_1(f,c)$ be an isotropic marginally trapped immersion of a surface into a Robertson-Walker spacetime. Then there are two possibilities.
\begin{itemize}
\item[(1)] the immersion takes the form
$$\phi = (t_0,\varphi),$$ 
where $t_0 \in I$ is constant and $\varphi$ defines a totally umbilical surface in $Q^3(c)$ with a mean curvature vector of constant length $|f'(t_0)|$;
\item[(2)] the function 
$$ \frac{\theta''\cos_c\theta+(\theta')^2c\sin_c\theta}{\theta''\sin_c\theta-(\theta')^2\cos_c\theta}, $$
where $\theta$ is defined in \eqref{eq:theta}, is constant, say $C_0$, and $\phi$ is locally congruent to an immersion of the form 
$$\phi = \left(\tau, \cos_c(\theta\circ\tau)\varphi+\sin_c(\theta\circ\tau)\nu\right),$$ 
where $\tau: S \to \R$ is a real function of class $C^2$ and $\varphi$ defines a totally umbilical surface in $Q^3(c)$ with Gauss map $\nu$ and with a mean curvature vector of length $|C_0|$.
\end{itemize}
\end{theorem}

%%%%%%%%%%%%%%%%%%%%%%%%%%%%%%%%%%%%%%%%%%%%%%%%%%%%%%%%%%%

\section{Marginally Trapped Surfaces with Constant Gaussian Curvature}
\label{sec:cgc}

In this section, we look for the marginally trapped surfaces with constant Gaussian curvature in the classifications from the previous sections. Note that, in theory, it suffices to find the constant Gaussian curvature surfaces in the theorems of Section \ref{sec:local_description}, but this seems to be a non-trivial task.

\subsection{Surfaces with Positive Relative Nullity} 

The following result follows from the equation of Gauss and the definition of positive relative nullity.

\begin{proposition} 
Every marginally trapped surface with positive relative nullity in a space form $Q^n_s(c)$, has constant Gaussian curvature $c$.
\end{proposition}

Therefore, all the surfaces listed in Theorem \ref{T:Positive:R}, Theorem \ref{T:Positive:S} and Theorem \ref{T:Positive:H} are examples of marginally trapped surfaces with constant Gaussian curvature.

\subsection{Surfaces in Light Cones}

Proposition \ref{P:Lightcones} showed that marginally trapped surfaces lying in the light cone of $Q^4_1(c)$, with $c\in\left\{-1,0,1 \right\}$, have constant Gaussian curvature $c$, while Proposition \ref{P:ExistenceLightcones} proved the existence of such surfaces. 

\subsection{Surfaces with Parallel Mean Curvature Vector Field} 

Marginally trapped surfaces with parallel mean curvature vector field in four-dimensional Lorentzian space forms are classified in Theorem \ref{T:ParallelMeanE}, Theorem \ref{T:Parallel:S} and Theorem \ref{T:Parallel:H}. The surfaces with constant Gaussian curvature are already mentioned in the formulation of these theorems and we can summarize the situation as follows.
\begin{itemize}
\item[(1)] Surfaces of types (1)--(4) in Theorem \ref{T:ParallelMeanE}, of types (2)--(4) in Theorem \ref{T:Parallel:S} and of types (2)--(4) in Theorem \ref{T:Parallel:H} are flat.
\item[(2)] Surfaces of type (1) in Theorem \ref{T:Parallel:S} have constant Gaussian curvature $K=1$.
\item[(3)] Surfaces of type (1) in Theorem \ref{T:Parallel:H} have constant Gaussian curvature $K=-1$.
\end{itemize}

\subsection{Boost Invariant Surfaces}

A boost invariant surface $S$ in Minkowski spacetime is locally congruent to a surface $\Sigma_\alpha$ with unit speed profile curve $\alpha=(\alpha_1, 0, \alpha_3, \alpha_4)$, where $\alpha_1$ is positive, and has Gaussian curvature
\begin{align*}
K=-\frac{{\alpha_1}''}{\alpha_1}.
\end{align*}
Theorem \ref{T:Boost} describes all boost invariant marginally trapped surfaces in $\mathbb{R}^4_1$ and the next proposition determines all such surfaces with constant Gaussian curvature.
 
\begin{proposition}\cite{Boost,Walleghem}
A boost invariant marginally trapped surface $S$ in Minkowski spacetime has constant Gaussian curvature if and only it is locally congruent to a surface $\Sigma_{\alpha}$ whose unit speed profile curve $\alpha=(\alpha_1, 0, \alpha_3, \alpha_4)$, with $\alpha_1 > 0$, is given by one of the following cases.
\begin{itemize}
\item[(1)]If $S$ is flat, then $\alpha$ is given by
\begin{align*}
\alpha_1(s)&=c_1 s+c_2,\\
\alpha_2(s)&=\frac{c_1 s+c_2}{\sqrt{c_1^2+1}}\left(\sin\left(\frac{\ln(c_1 s+c_2)}{c_1}+\xi_0 \right) +c_1\cos\left(\frac{\ln(c_1 s+c_2)}{c_1}+\xi_0 \right)\right),\\
\alpha_3(s)&=\frac{c_1 s+c_2}{\sqrt{c_1^2+1}}\left(\cos\left(\frac{\ln(c_1 s+c_2)}{c_1}+\xi_0 \right)-c_1\sin\left(\frac{\ln(c_1 s+c_2)}{c_1}+\xi_0 \right) \right),
\end{align*}
with $c_1, c_2,\xi_0\in\mathbb{R}$.
\item[(2)]If $S$ has constant Gaussian curvature $K>0$, then $\alpha$ is given by
\begin{align*}
\alpha_1(s)&=c_1\cos(\sqrt{K}s+c_2),\\
\alpha_3(s)&= \int_{s_0}^s \sqrt{1+c_1^2K\sin^2(\sqrt{K}t+c_2)} \cos\xi(t) \, \de t,\\
\alpha_4(s) &= \int_{s_0}^s \sqrt{1+c_1^2K\sin^2(\sqrt{K}t+c_2)}\sin\xi(t) \, \de t, 
\end{align*}
where $c_1, c_2\in\mathbb{R},s_0\in I$ and the function $\xi$ is given by
\begin{align*}
\xi(t)=\int_{s_0}^t \frac{c_1^2K\cos^2(\sqrt{K}u+c_2)-2}{c_1\cos(\sqrt{K}u+c_2)(1+c_1^2K\sin(\sqrt{K}u+c_2))} \, \de u.
\end{align*}
\item[(3)]
If $S$ has constant Gaussian curvature $K<0$, then $\alpha$ is given by
\begin{align*}
\alpha_1(s)&=c_1\exp(\sqrt{-K}s)+c_2\exp(-\sqrt{-K}s),\\ 
\alpha_3(s)&= \int_{s_0}^s \sqrt{1-K\left(c_1\exp(\sqrt{-K}t)+c_2\exp(-\sqrt{K}t)\right)^2} \cos\xi(t) \, \de t,\\
\alpha_4(s) &= \int_{s_0}^s \sqrt{1-K\left(c_1\exp(\sqrt{-K}t)+c_2\exp(-\sqrt{K}t)\right)^2}\sin\xi(t) \, \de t,
\end{align*}
where $c_1, c_2\in\mathbb{R},s_0\in I$ and the function $\xi$ is given by
\begin{align*}
\xi(t)=\int_{s_0}^t\frac{1-c_1^2K\exp(2\sqrt{-K}u)-c_2^2K\exp(-2\sqrt{-K}u)}{1-K\left(c_1\exp(\sqrt{-K}u)+c_2\exp(-\sqrt{K}u)\right)^2} \, \de u.
\end{align*}
\end{itemize}
\end{proposition}

\subsection{Rotation Invariant Surfaces}

A rotation invariant surface $S$ in Minkowski spacetime is locally congruent to a surface $\Sigma_\alpha$ with unit profile curve $\alpha=(\alpha_1, \alpha_2, \alpha_3, 0)$, where $\alpha_3$ is positive, and has Gaussian curvature
\begin{align*}
K=-\frac{{\alpha_3}''}{\alpha_3}.
\end{align*}
Theorem \ref{T:Rotation} classifies all rotation invariant marginally trapped surfaces in $\mathbb{R}^4_1$ and the next proposition determines all such surfaces with constant Gaussian curvature.

\begin{proposition}\cite{Rotation}
A rotation invariant marginally trapped surface $S$ in Minkowski spacetime has constant Gaussian curvature if and only it is locally congruent to a surface $\Sigma_{\alpha}$ whose unit speed profile curve $\alpha=(\alpha_1, \alpha_2, \alpha_3, 0)$, with $\alpha_3 > 0$, is given by one of the following cases.
\begin{itemize}
\item[(1)]If $S$ is flat, then $\alpha$ is given by one of the following curves:
\begin{itemize}
\item[(i)] a curve of type (1) in Theorem \ref{T:Rotation};
\item[(ii)]
a curve of type (2) in Theorem \ref{T:Rotation}, described by one of the following two cases:
\begin{itemize}
\item[(a)]the coordinate functions are given by
\begin{align*}
\alpha_1(s)&=\frac{1}{2}\left(\frac{1-c_1}{1+c_1}(c_1s+c_2)^{\frac{c_1+1}{c_1}}+\frac{1+c_1}{1-c_1}(c_1s+c_2)^{\frac{c_1-1}{c_1}} \right),\\
\alpha_2(s)&=\frac{1}{2}\left(\frac{1-c_1}{1+c_1}(c_1s+c_2)^{\frac{c_1+1}{c_1}}-\frac{1+c_1}{1-c_1}(c_1s+c_2)^{\frac{c_1-1}{c_1}} \right),\\
\alpha_3(s)&=c_1s+c_2,
\end{align*}
where $c_1,c_2\in\mathbb{R}$, with $|c_1|\notin\left\{0,1\right\}$,
\vskip.08in
\item[(b)]$\displaystyle{\alpha(s)=\left(c_1\cosh\left(\frac{s}{c_1} \right),  c_1\sinh\left(\frac{s}{c_1} \right),c_1,0 \right)}$, with $c_1\in\mathbb{R}_0^{+}$;
\end{itemize}
\end{itemize}
\item[(2)]If $S$ has constant Gaussian curvature $K>0$, then $\alpha$ is of type (2) in Theorem \ref{T:Rotation}, with
\begin{align*}
\alpha_1(s)&=\left(\sinh\xi(s)+c_1\sin(\sqrt{K}s+c_2)\cosh\xi(s) \right),\\
\alpha_2(s)&=\left(\cosh\xi(s)+c_1\sin(\sqrt{K}s+c_2)\sinh\xi(s) \right),\\
\alpha_3(s)&=c_1\cos(\sqrt{K}s+c_2),
\end{align*} 
where $c_1, c_2\in\mathbb{R}, c_1\neq0$ and the function $\xi$ is given by
\begin{align*}
\xi(s)=\frac{1}{c_1}\ln\left|\frac{1+\sin(\sqrt{K}s+c_2)}{1-\sin(\sqrt{K}s+c_2)} \right|,
\end{align*}
with $\xi_0\in\mathbb{R}$.
\item[(3)]If $S$ has constant Gaussian curvature $K<0$, then $\alpha$ is of type (2) in Theorem \ref{T:Rotation}, with
\begin{align*}
\alpha_1(s)&=\left(\sinh\xi(s)-\sqrt{-K}(c_1-c_2)\exp\xi(s)\cosh\xi(s) \right),\\
\alpha_2(s)&=\left(\cosh\xi(s)-\sqrt{-K}(c_1-c_2)\exp\xi(s)\sinh\xi(s) \right),\\
\alpha_3(s)&=c_1\exp(\sqrt{-K}s)+c_2\exp(-\sqrt{-K}s),
\end{align*}
where $ c_1, c_2\in\mathbb{R}$ with $c_1^2+c_2^2>0$ and the fuction $\xi$ is given by
\begin{alignat*}{3}
\xi(s)&=\frac{1}{\sqrt{-c_1c_2K}}\arctan\left(\frac{c_1\exp(\sqrt{-K}s)}{\sqrt{c_1c_2}} \right)+\xi_0\qquad&&\mbox{ if }c_1c_2>0,\\
\xi(s)&=\frac{1}{\sqrt{2c_1c_2K}}\ln\left|\frac{2c_1\exp(\sqrt{-K}s)-2\sqrt{-c_1c_2}}{\sqrt{-K}s+2\sqrt{-c_1c_2}} \right|+\xi_0&&\mbox{ if }c_1c_2<0,\\
\xi(s)&=-\frac{1}{c_1\sqrt{-K}\exp(\sqrt{-K}s)}+\xi_0&&\mbox{ if }c_2=0,\\
\xi(s)&=\frac{\exp(\sqrt{-K}s)}{c_2\sqrt{-K}}+\xi_0&&\mbox{ if }c_1=0,
\end{alignat*}
with $\xi_0\in\mathbb{R}$.
\end{itemize}
\end{proposition}

\subsection{Screw Invariant Surfaces}

A screw invariant surface $S$ in Minkowski spacetime is locally congruent to a surface $\Sigma_\alpha$ with unit speed profile curve $\alpha=(\alpha_k, \alpha_l, 0, \alpha_4)$, where $\alpha_k$ is positive, and has a Gaussian curvature described by 
\begin{align*}
K=-\frac{{\alpha_k}''}{\alpha_k}.
\end{align*}
Theorem \ref{T:Screw} classifies all screw invariant marginally trapped surfaces in $\mathbb{R}^4_1$ and the next proposition determines all such surfaces with constant Gaussian curvature.

\begin{proposition}\cite{Screw}
A screw invariant marginally trapped surface $S$ in Minkowski spacetime has constant Gaussian curvature if and only it is locally congruent to a surface $\Sigma_{\alpha}$ whose unit speed profile curve $\alpha=(\alpha_k, \alpha_l, 0, \alpha_4)$, with $\alpha_k$ positive, is given by one of the following cases.
\begin{itemize}
\item[(1)] If $S$ is flat, then $\alpha$ is a curve of type (1) in Theorem \ref{T:Screw}.
\item[(2)] If $S$ has constant Gaussian curvature $K>0$, then $\alpha$ is a curve of type (2) in Theorem \ref{T:Screw} with 
\begin{align*}
\rho(s)=2c_1^2K\cos\left(2\sqrt{K}s+c_2 \right),
\end{align*}
with $c_1,c_2\in\mathbb{R}$.
\item[(3)] If $S$ has constant Gaussian curvature $K<0$, then $\alpha$ is a curve of type (2) in Theorem \ref{T:Screw} with
\begin{align*}
\rho(s)=\frac{2\sqrt{-K}\left(c_1^2\exp(2\sqrt{-K}s)+c_2^2\exp(-2\sqrt{-K}s) \right)}{c_1^2\exp(2\sqrt{-K}s)-c_2^2\exp(-2\sqrt{-K}s) },
\end{align*}
with $c_1,c_2\in\mathbb{R}$.
\end{itemize}
\end{proposition}

\subsection{Isotropic surfaces} 

Theorem \ref{T:isotropy} implies that an isotropic marginally trapped surface in a Lorentzian space form $Q^4_1(c)$ has constant Gaussian curvature $c$. All isotropic marginally trapped surfaces in these spacetimes therefore provide examples of constant Gaussian curvature marginally trapped surfaces. 

%%%%%%%%%%%%%%%%%%%%%%%%%%%%%%%%%%%%%%%%%%%%%%%%%%%%%%%%%%%

\section{Conclusions and Open Questions}
\label{sec:openproblems}

While trapped and marginally trapped surfaces are concepts from physics, they have very natural geometric definitions. In particular, a marginally trapped surface is a Riemannian surface in a spacetime whose mean curvature vector field, one of the most important invariants in submanifold geometry, is lightlike at every point. It is hence no surprise that this family of surfaces and their generalizations were studied intensively from a purely geometric point of view. In this paper we gave an overview of this study when the ambient space is Minkowski space, de Sitter space, anti-de Sitter space or a Robertson-Walker spacetime. Most results are classification theorems under additional geometric conditions. The local descriptions in Section \ref{sec:local_description} provide in some sense a complete classification without additional assumptions, but it is not always easy to find surfaces with particular properties from this general description, such as constant Gaussian curvature surfaces, see Section \ref{sec:cgc}. 

We finish the paper with some open questions regarding marginally trapped surfaces in spacetimes, which arise naturally form the current overview article.

\begin{itemize}
\item[(1)] What are the marginally trapped surfaces with parallel mean curvature vector in a Robertson-Walker spacetime?
\item[(2)] For any one-parameter group of isometries of de Sitter, anti-de Sitter or a Robertson-Walker spacetime: what are the invariant marginally trapped surfaces?
\item[(3)] What are the marginally trapped surfaces with constant Gaussian curvature in Minkowski, de Sitter, anti-de Sitter or a Robertson-Walker spacetime? In particular, what are the flat marginally trapped surfaces? 
\item[(4)] What are the marginally trapped surfaces satisfying any of the additional conditions appearing in this paper in other (four-dimensional) Lorentzian manifolds, such as Kerr spacetime and Schwarzschild spacetime?
\end{itemize}

Remark that one could in principle start from the general descriptions given in Section \ref{sec:local_description} to tackle questions (1)--(3).
\vspace{6PT}

%%%%%%%%%%%%%%%%%%%%%%%%%%%%%%%%%%%%%%%%%%%%%%%%%%%%%%%%%%%

\bibliography{citations}
\bibliographystyle{plain}

\end{document}